\input amstex
\documentstyle{amsppt}
%----------------------------------------------------------------
% Title:     Hadamard matrices in $\{0,1\}$ presentation and
%            an algorithm for generating  them.
% Author:    Ruslan Sharipov
% Comments:  AmSTeX, 12 pages, amsppt style
% MSC-class: 05B20, 11D04, 11D09, 15B10, 15B34, 15B36, 65-04
%----------------------------------------------------------------
%           Replacement for output macro definition
%
\catcode`@=11
\redefine\output@{%
  \def\break{\penalty-\@M}\let\par\endgraf
  \ifodd\pageno\global\hoffset=105pt\else\global\hoffset=8pt\fi  
  \shipout\vbox{%
    \ifplain@
      \let\makeheadline\relax \let\makefootline\relax
    \else
      \iffirstpage@ \global\firstpage@false
        \let\rightheadline\frheadline
        \let\leftheadline\flheadline
      \else
        \ifrunheads@ %\let\makefootline\relax
        \else \let\makeheadline\relax
        \fi
      \fi
    \fi
    \makeheadline \pagebody \makefootline}%
  \advancepageno \ifnum\outputpenalty>-\@MM\else\dosupereject\fi
}
\def\Beta{\mathchar"0\hexnumber@\rmfam 42}
\redefine\mm@{2010} % Math Subject Classification year
\catcode`\@=\active
%----------------------------------------------------------------
\nopagenumbers
\chardef\textvolna='176

\chardef\bigalpha='013
\def\negskp{\hskip -2pt}

\chardef\degree="5E
\def\compos{\,\raise 1pt\hbox{$\sssize\circ$} \,}
%\def\id{\operatorname{id}}

%\font\eightrm=cmr8
%\def\LT{\operatorname{\text{\eightrm LT}}}
%\def\LM{\operatorname{\text{\eightrm LM}}}
%\def\LC{\operatorname{\text{\eightrm LC}}}
%\accentedsymbol\hatgamma{\kern 2pt\hat{\kern -2pt\gamma}}
%\accentedsymbol\checkgamma{\kern 2.5pt\check{\kern -2.5pt\gamma}}
\def\blue#1{#1}

\gdef\darkred#1{#1}
\gdef\darkblue#1{#1}
\catcode`#=11\def\diez{#}\catcode`#=6
\catcode`&=11\catcode`&=4
\catcode`_=11\def\podcherkivanie{_}\catcode`_=8
\catcode`\^=11\catcode`\^=7
\catcode`~=11\catcode`~=\active
\catcode`\%=11\def\procent{%}\catcode`\%=14
\def\mycite#1{\cite{\blue{#1}}\immediate\special{ps:
     ShrHPSdict begin /ShrBORDERthickness 0 def}}
\def\myciterange#1#2#3#4{\cite{\blue{#2#3#4}}\immediate\special{ps:
     ShrHPSdict begin /ShrBORDERthickness 0 def}}
\def\mytag#1{%
    \tag#1}
\def\mythetag#1{\thetag{\blue{#1}}\immediate\special{ps:
     ShrHPSdict begin /ShrBORDERthickness 0 def}}
\def\myrefno#1{\no#1}
\def\myhref#1#2{\blue{#2}\immediate\special{ps:
     ShrHPSdict begin /ShrBORDERthickness 0 def}}
\def\myEarXivlink{\myhref{http://arXiv.org}{http:/\negskp/arXiv.org}}
\def\myGeoCities{\myhref{http://www.geocities.com}{GeoCities}}
\def\mytheorem#1{\csname proclaim\endcsname{Theorem #1}}
\def\mytheoremwithtitle#1#2{\csname proclaim\endcsname{Theorem #1#2}}
\def\mythetheorem#1{\blue{#1}\immediate\special{ps:
     ShrHPSdict begin /ShrBORDERthickness 0 def}}
\def\mylemma#1{\csname proclaim\endcsname{Lemma #1}}
\def\mylemmawithtitle#1#2{\csname proclaim\endcsname{Lemma #1#2}}
\def\mythelemma#1{\blue{#1}\immediate\special{ps:
     ShrHPSdict begin /ShrBORDERthickness 0 def}}
\def\mycorollary#1{\csname proclaim\endcsname{Corollary #1}}
\def\mythecorollary#1{\blue{#1}\immediate\special{ps:
     ShrHPSdict begin /ShrBORDERthickness 0 def}}
\def\mydefinition#1{\definition{Definition #1}}
\def\mythedefinition#1{\blue{#1}\immediate\special{ps:
     ShrHPSdict begin /ShrBORDERthickness 0 def}}
\def\myconjecture#1{\csname proclaim\endcsname{Conjecture #1}}
\def\myconjecturewithtitle#1#2{\csname proclaim\endcsname{Conjecture #1#2}}
\def\mytheconjecture#1{\blue{#1}\immediate\special{ps:
     ShrHPSdict begin /ShrBORDERthickness 0 def}}
\def\myproblem#1{\csname proclaim\endcsname{Problem #1}}
\def\myproblemwithtitle#1#2{\csname proclaim\endcsname{Problem #1#2}}
\def\mytheproblem#1{\blue{#1}\immediate\special{ps:
     ShrHPSdict begin /ShrBORDERthickness 0 def}}
\def\mytable#1{Table #1}
\def\mythetable#1{\blue{#1}\immediate\special{ps:
     ShrHPSdict begin /ShrBORDERthickness 0 def}}
\def\myanchortext#1#2{#2}
\def\mytheanchortext#1#2{\blue{#2}\immediate\special{ps:
     ShrHPSdict begin /ShrBORDERthickness 0 def}}
%----------------------------------------------------------------
% Cyrillic fonts definition
\font\eightcyr=wncyr8
%\font\tencyr=wncyr10%----------------------------------------------------------------
\pagewidth{360pt}
\pageheight{606pt}
\topmatter
\title
Hadamard matrices in \{0,1\} presentation and an
algorithm for generating  them.
\endtitle
\rightheadtext{Hadamard matrices in \{0,1\} presentation}
\author
Ruslan Sharipov
\endauthor
\address Bashkir State University, 32 Zaki Validi street, 450074 Ufa, Russia
\endaddress
\email
\myhref{mailto:r-sharipov\@mail.ru}{r-sharipov\@mail.ru}
\endemail
\abstract
Hadamard matrices are square $n\times n$ matrices whose entries are ones and minus 
ones and whose rows are orthogonal to each other with respect to the
standard scalar product in $\Bbb R^n$. Each Hadamard matrix can be transformed to
a matrix whose entries are zeros and ones. This presentation of Hadamard matrices 
is investigated in the paper and based on it an algorithm for generating them is 
designed. 
\endabstract
\subjclassyear{2010}
\subjclass 05B20, 11D04, 11D09, 15B10, 15B34, 15B36, 65-04\endsubjclass
\keywords Hadamard matrices
\endkeywords
\endtopmatter
\loadbold
%\loadeufb
\TagsOnRight
\document

\head
1. Introduction.
\endhead
     Hadamard matrices are known for $n=1$, for $n=2$, and for $n=4\,q$, where 
$q\in\Bbb N$ and $\Bbb N$ is the set of positive integers. However it is not yet 
known if they do exist for all $q\in\Bbb N$ (see \mycite{1}). Hadamard matrices 
are associated with Hadamard's maximal determinant problem (see \mycite{2} and
\mycite{3}) and solve this problem for $n=4\,q$, where $q\in\Bbb N$. For the 
general case $n\in\Bbb N$ Hadamard's maximal determinant problem is yet unsolved. 
Its simplified version is suggested in \mycite{4}.\par
     Let $H$ be an $n\times n$ Hadamard matrix. By definition its rows considered
as vectors in $\Bbb R^n$ are orthogonal to each other with respect to the standard 
scalar product in $\Bbb R^n$. The same is valid for its columns. The proof is 
elementary. Indeed, since $|\bold r_i|^2=n$ for each row $\bold r_i$ treated as a 
vector, the orthogonality of rows implies
$$
\hskip -2em
\sum^n_{k=1}H_{ik}\,H_{jk}=n\,\delta_{ij},
\mytag{1.1}
$$
where $\delta_{ij}$ is the Kronecker delta (see \S\,23 in
Chapter~\uppercase\expandafter{\romannumeral 1} of \mycite{5}). The equality 
\mythetag{1.1} means $H\,H^{\top}=n\,I$, where $H^{\top}$ is the transpose of $H$ and 
$I$ is the identity matrix. Then $H\,(H^{\top}/n)=I$ and $H^{-1}=H^{\top}/n$, which 
yields the equalities $(H^{\top}/n)\,H=I$ and $H^{\top}\,H=n\,I$. The latter one is 
written as
$$
\hskip -2em
\sum^n_{k=1}H_{ki}\,H_{kj}=n\,\delta_{ij}.
\mytag{1.2}
$$
The equality \mythetag{1.2} is exactly the orthogonality condition for the columns 
of $H$.\par
    For each particular $n\in\Bbb N$ the set of $n\times n$ Hadamard matrices is 
invariant under the following transformations (see \mycite{6}):
\roster
\item"1)" permutation of rows/columns; 
\item"2)" multiplication of any row/column by $-1$.
\endroster 
Using these transformations each Hadamard matrix can be brought to a form where its 
first row and its first column both are filled with ones only (see \mycite{6}):
$$
\hskip -2em
H=\Vmatrix 1 & 1 & 1 &\hdots &1\\
1 & \vtop{\hsize 55pt\vskip -6pt
\leftline{\vphantom{a}\kern -4pt
\boxed{\vtop to 38pt{\hsize 50pt
\vskip 15pt\centerline{$\{-1,1\}$}
\vss}}}
\vskip -40pt\vss
}\kern -52pt\\
1\\
\vdots\\
1
\endVmatrix.
\mytag{1.3}
$$
Almost all Hadamard matrices in Sloan's library \mycite{7} are presented in the form 
\mythetag{1.3}.\par
    By writing $\{-1,1\}$ in \mythetag{1.3} we indicate that the lower right block of
the matrix is built by ones and minus ones. The next trick is to subtract the first 
row from each of the other rows in \mythetag{1.3}. As a result we get the matrix
$$
\hskip -2em
M=\Vmatrix 1 & 1 & 1 &\hdots &1\\
0 & \vtop{\hsize 55pt\vskip -6pt
\leftline{\vphantom{a}\kern -4pt
\boxed{\vtop to 38pt{\hsize 50pt
\vskip 15pt\centerline{$\{-2,0\}$}
\vss}}}
\vskip -40pt
}\kern -52pt\\
0\\
\vdots\\
0
\endVmatrix.
\mytag{1.4}
$$
Let's denote through $\tilde H$ the lower right block of the matrix 
\mythetag{1.4} divided by $-2$:
$$
\hskip -2em
\tilde H=\frac{1}{-2}\cdot
\lower 15pt\hbox{\boxed{\vbox to 38pt{\hsize 50pt
\vskip 15pt\centerline{$\{-2,0\}$}
\vss}}}
=\lower 15pt\hbox{\boxed{\vbox to 38pt{\hsize 50pt
\vskip 15pt\centerline{$\{0,1\}$}
\vss}}}\ .
\mytag{1.5}
$$
The transformation $H\longrightarrow\tilde H$ given by \mythetag{1.3}, \mythetag{1.4}, 
\mythetag{1.5} is well-known (see \mycite{2}). 
\mydefinition{1.1} The matrix $\tilde H$ produced from a Hadamard matrix $H$ 
of the form \mythetag{1.3} according to \mythetag{1.4} and \mythetag{1.5} is
called the $\{0,1\}$ presentation of the matrix $H$.
\enddefinition
    It is obvious that the transformation $H\longrightarrow\tilde H$ is invertible, 
i\.\,e\. each matrix $H$ of the form \mythetag{1.3} is associated with a 
unique matrix $\tilde H$ of the form \mythetag{1.5} and vice versa. Our goal
in this paper is to study $\{0,1\}$ presentation of Hadamard matrices and 
to design an algorithm for generating them in this presentation.\par
\head
2. Gram matrices. 
\endhead
\mydefinition{2.1} For any ordered list of vectors $\bold e_1,\,\ldots,\,
\bold e_s$ in a Euclidean space $\Bbb E$ their Gram matrix is the following 
matrix composed by their mutual scalar products (see \mycite{8} or \S\,1 in 
Chapter~\uppercase\expandafter{\romannumeral 5} of \mycite{9}):
$$
G=\Vmatrix (\bold e_1,\bold e_1) & \hdots & (\bold e_1,\bold e_s)\\
\vdots & \ddots & \vdots\\
(\bold e_s,\bold e_1) & \hdots & (\bold e_s,\bold e_s)
\endVmatrix. 
$$
\enddefinition
Let $\tilde\bold r_1,\,\ldots,\,\tilde\bold r_{n-1}$ be rows of the
matrix $\tilde H$ in \mythetag{1.5} considered as vectors in $\Bbb R^{n-1}$. 
If we enumerate the entries of  $H$ in \mythetag{1.3} starting from zero, 
then we can write 
$$
\hskip -2em
\tilde\bold r_i=f(\bold r_i), \ \ i=1,\,\ldots,\,n-1,
\mytag{2.1}
$$
where $\bold r_0,\,\bold r_1,\,\ldots,\,\bold r_{n-1}$ are the rows of 
the matrix \mythetag{1.3} and $f\!:H\longrightarrow\tilde H$ is the mapping given 
by \mythetag{1.3}, \mythetag{1.4}, and \mythetag{1.5}. From \mythetag{2.1} 
we derive
$$
\hskip -2em
(\tilde\bold r_i,\tilde\bold r_j)=\Bigl(\frac{\bold r_i-\bold r_0}{-2},
\frac{\bold r_j-\bold r_0}{-2}\Bigr).
\mytag{2.2}
$$
The equality \mythetag{2.2} holds since the initial entry in almost all rows
of the matrix \mythetag{1.4} is zero. Expanding the right hand side of 
\mythetag{2.2}, we get
$$
\hskip -2em
(\tilde\bold r_i,\tilde\bold r_j)=\frac{(\bold r_i,\bold r_j)}{4}
-\frac{(\bold r_0,\bold r_i)}{4}-\frac{(\bold r_0,\bold r_j)}{4}+
\frac{(\bold r_0,\bold r_0)}{4}.
\mytag{2.3}
$$
Since $H$ is a Hadamard matrix, taking into account the renumeration of the 
entries of $H$, from \mythetag{1.1} we derive the following four equalities:
$$
\xalignat 2
&\hskip -2em
(\bold r_i,\bold r_j)=n\,\delta_{ij},
&&(\bold r_i,\bold r_0)=0,\\
\vspace{-1.5ex}
\mytag{2.4}\\
\vspace{-1.5ex}
&\hskip -2em
(\bold r_j,\bold r_0)=0,
&&(\bold r_0,\bold r_0)=n.
\endxalignat
$$
Applying \mythetag{2.4} to \mythetag{2.3} yields 
$$
\hskip -2em
(\tilde\bold r_i,\tilde\bold r_j)=\frac{n\,(\delta_{ij}+1)}{4}.
\mytag{2.5}
$$
The equality \mythetag{2.5} means that we have proved the following
theorem.
\mytheorem{2.1} For any $m\times m$ Hadamard matrix in $\{0,1\}$ presentation
with $m>1$ its size $m=4\,q-1$, where $q\in\Bbb N$, and the Gram matrix associated
with the rows of this Hadamard matrix looks like
$$
\hskip -2em
G=\Vmatrix 
b & a & \hdots & a\\
a & b & \hdots & a\\
\vdots & \vdots & \ddots &\vdots\\
a & a & \hdots & b
\endVmatrix,
\mytag{2.6}
$$
where $a=(m+1)/4=q$ and $b=(m+1)/2=2\,q$.
\endproclaim
    The calculations \mythetag{2.2}, \mythetag{2.3}, \mythetag{2.4}, and \mythetag{2.5}
are invertible. Therefore Theorem~\mythetheorem{2.1} can be strengthened as follows.
\mytheorem{2.2} A square $m\times m$ matrix with $m=4\,q-1$ whose entries are zeros and
ones is a Hadamard matrix in $\{0,1\}$ presentation if an a only if the Gram matrix 
associated with its rows is of the form \mythetag{2.6}, where $a=(m+1)/4=q$ and $b=(m+1)/2=2\,q$.
\endproclaim
\demo{Proof} The necessity part in the statement of Theorem~\mythetheorem{2.2} is proved
by Theorem~\mythetheorem{2.1}. Let's prove the sufficiency.\par
    Going backward from \mythetag{1.5} to \mythetag{1.4}, we insert the initial column of 
zeros to $\tilde H$. Since $|\tilde\bold r_i|^2=b=(m+1)/2=2\,q=n/2$ in \mythetag{2.6}, 
upon this step we get an $m\times n$ matrix with equal number of zeros and ones in each
row. Then we multiply this matrix by $-2$, insert the initial row composed by ones, and add 
this initial row to each of the other rows. As a result we get a matrix of the form
\mythetag{1.3}. Each of its rows, except for the initial one, comprises equal number of 
ones and minus ones. This means that the equalities $(\bold r_i,\bold r_0)=0$ and
$(\bold r_j,\bold r_0)=0$ in \mythetag{2.4} are fulfilled. The equality $(\bold r_0,
\bold r_0)=n$ in \mythetag{2.4} holds since the initial row $\tilde r_0$ in \mythetag{1.3} 
is composed by ones only. The equality $(\bold r_i,\bold r_j)=n\,\delta_{ij}$ in
\mythetag{2.4} is derived from \mythetag{2.5} using \mythetag{2.3}, while \mythetag{2.5}
itself follows from \mythetag{2.6} since $a=(m+1)/4=q=n/4$ and $b=(m+1)/2=2\,q=n/2$.
The whole set of the equalities \mythetag{2.4} is equivalent to \mythetag{1.1} upon passing
to the standard enumeration of the entries of $H$, thus proving that the matrix $H$ in
\mythetag{1.3} produced backward from \mythetag{1.5} through \mythetag{1.4} is a regular
Hadamard matrix. Theorem~\mythetheorem{2.2} is proved.\qed\enddemo
     Now let's consider the columns of the matrix $\tilde H$ in \mythetag{1.5}. We denote
them through $\tilde\bold c_1,\,\ldots,\,\tilde\bold c_{n-1}$. If we enumerate the entries
of the matrix $H$ in \mythetag{1.3} starting from zero, then for $\tilde\bold c_1,\,\ldots,
\,\tilde\bold c_{n-1}$ we can write 
$$
\hskip -2em
\tilde\bold c_i=f(\bold c_i), \ \ i=1,\,\ldots,\,n-1,
\mytag{2.7}
$$
where $\bold c_0,\,\bold c_1,\,\ldots,\,\bold c_{n-1}$ are the columns of the matrix 
\mythetag{1.3} and $f\!:H\longrightarrow\tilde H$ is the mapping given by \mythetag{1.3}, \mythetag{1.4}, and \mythetag{1.5}. From \mythetag{2.7} we derive
$$
(\tilde\bold c_i,\tilde\bold c_j)=\sum^{n-1}_{k=1}\frac{(H_{ki}-H_{0i})}{-2}\,
\frac{(H_{kj}-H_{0j})}{-2}=\sum^{n-1}_{k=1}\frac{(H_{ki}-H_{0i})\,(H_{kj}-H_{0j})}
{4}\,.
$$
The right hand side of this equality can be transformed as 
$$
\sum^{n-1}_{k=1}\frac{(H_{ki}-H_{0i})\,(H_{kj}-H_{0j})}{4}
=\sum^{n-1}_{k=0}\frac{(H_{ki}-H_{0i})\,(H_{kj}-H_{0j})}{4}
$$
since the term with $k=0$ in the sum does vanish. As a result we get   
$$
\hskip -2em
(\tilde\bold c_i,\tilde\bold c_j)=
\sum^{n-1}_{k=0}\frac{(H_{ki}-H_{0i})\,(H_{kj}-H_{0j})}{4}.
\mytag{2.8}
$$
But $H_{0i}=H_{k0}=1$ and $H_{0j}=H_{k0}=1$ due to \mythetag{1.3}. Therefore from
\mythetag{2.8} we get 
$$
\hskip -2em
(\tilde\bold c_i,\tilde\bold c_j)=
\sum^{n-1}_{k=0}\frac{(H_{ki}-H_{k0})\,(H_{kj}-H_{k0})}{4}.
\mytag{2.9}
$$
Expanding the right hand side of \mythetag{2.9}, we write
$$
\hskip -2em
\gathered
(\tilde\bold c_i,\tilde\bold c_j)=\sum^{n-1}_{k=0}\frac{H_{ki}\,H_{kj}}{4}
-\sum^{n-1}_{k=0}\frac{H_{ki}\,H_{k0}}{4}\,-\\
-\sum^{n-1}_{k=0}\frac{H_{kj}\,H_{k0}}{4}+\sum^{n-1}_{k=0}\frac{H_{k0}\,H_{k0}}{4}.
\endgathered
\mytag{2.10}
$$
Each sum in \mythetag{2.10} is expressed through the scalar product of some definite 
pair of columns of the matrix \mythetag{1.3}. Indeed we have
$$
\hskip -2em
\gathered
(\tilde\bold c_i,\tilde\bold c_j)=\frac{(\bold c_i,\bold c_j)}{4}
-\frac{(\bold c_i,\bold c_0)}{4}-\frac{(\bold c_j,\bold c_0)}{4}
+\frac{(\bold c_0,\bold c_0)}{4}.
\endgathered
\mytag{2.11}
$$
The equality \mythetag{2.11} is similar to \mythetag{2.3}. Since $H$ is a Hadamard 
matrix, taking into account the renumeration of the entries of $H$, from \mythetag{1.2} 
we derive
$$
\xalignat 2
&\hskip -2em
(\bold c_i,\bold c_j)=n\,\delta_{ij},
&&(\bold c_i,\bold c_0)=0,\\
\vspace{-1.5ex}
\mytag{2.12}\\
\vspace{-1.5ex}
&\hskip -2em
(\bold c_j,\bold c_0)=0,
&&(\bold c_0,\bold c_0)=n.
\endxalignat
$$
Applying \mythetag{2.12} to \mythetag{2.11}, we derive an equality which is similar
to \mythetag{2.5}:
$$
\hskip -2em
(\tilde\bold c_i,\tilde\bold c_j)=\frac{n\,(\delta_{ij}+1)}{4}.
\mytag{2.13}
$$
The equality \mythetag{2.13} means that we have proved a theorem similar to
Theorem~\mythetheorem{2.1}. 
\mytheorem{2.3} For any $m\times m$ Hadamard matrix in $\{0,1\}$ presentation
with $m>1$ its size $m=4\,q-1$, where $q\in\Bbb N$, and the Gram matrix associated
with the columns of this Hadamard matrix looks like
$$
\hskip -2em
G=\Vmatrix 
b & a & \hdots & a\\
a & b & \hdots & a\\
\vdots & \vdots & \ddots &\vdots\\
a & a & \hdots & b
\endVmatrix,
\mytag{2.14}
$$
where $a=(m+1)/4=q$ and $b=(m+1)/2=2\,q$.
\endproclaim
    Note that the Gram matrices in \mythetag{2.14} and \mythetag{2.6} do coincide
though their entries are defined differently. A theorem similar to 
Theorem~\mythetheorem{2.2} is formulated as follows. 
\mytheorem{2.4} A square $m\times m$ matrix with $m=4\,q-1$ whose entries are zeros and
ones is a Hadamard matrix in $\{0,1\}$ presentation if an a only if the Gram matrix 
associated with its columns is of the form \mythetag{2.6}, where $a=(m+1)/4=q$ and $b=(m+1)/2=2\,q$.
\endproclaim
\demo{Proof} The necessity part in the statement of Theorem~\mythetheorem{2.4} is proved
by Theorem~\mythetheorem{2.3}. Let's prove the sufficiency.\par
     When producing the matrix \mythetag{1.3} backward from the matrix \mythetag{1.5}
each one in the matrix \mythetag{1.5} becomes minus one in the matrix \mythetag{1.3}
and each zero in the matrix \mythetag{1.5} becomes one in the matrix \mythetag{1.3}.
Extra ones in the initial row and in the initial column of the matrix \mythetag{1.5}
are inserted independently. Therefore the equality $|\tilde\bold c_i|^2=b=(m+1)/2
=2\,q=n/2$ for the diagonal entries in \mythetag{2.14} means that the number of ones
is equal to the number of minus ones in each column of the matrix \mythetag{1.5}, except
for the initial column $\bold c_0$. This yields the equalities $(\bold c_i,\bold c_0)=0$
and $(\bold c_j,\bold c_0)=0$ in \mythetag{2.12}. The equality $(\bold c_0,\bold c_0)=n$
in \mythetag{2.12} holds since the initial column of the matrix \mythetag{1.3} is composed
by ones only. Then the equality $(\bold c_i,\bold c_j)=n\,\delta_{ij}$ in \mythetag{2.12}
is derived from \mythetag{2.13} using \mythetag{2.11}, while \mythetag{2.13} itself follows 
from \mythetag{2.14} since $a=(m+1)/4=q=n/4$ and $b=(m+1)/2=2\,q=n/2$. The whole set of the equalities \mythetag{2.12} is equivalent to \mythetag{1.2} upon passing to the standard enumeration of the entries of $H$, thus proving that the matrix $H$ in \mythetag{1.3} 
produced backward from \mythetag{1.5} through \mythetag{1.4} is a regular
Hadamard matrix.\qed\enddemo
\head
3. Partitions and groupings\\
in rows of $\{0,1\}$ Hadamard matrices. 
\endhead
     Due to Theorem~\mythetag{2.2} and Theorem~\mythetag{2.4} the whole set of Hadamard
matrices in $\{0,1\}$ presentation is invariant under permutation of rows and columns 
of the matrices. Any two matrices produced from each other by means of these 
transformations are called equivalent. Below they are threated as inessentially 
different.\par
     Let $H$ be some particular $15\times 15$ Hadamard matrix in $\{0,1\}$ presentation. 
We choose the following one as an example:
$$
H=\Vmatrix
1&1&1&1&1&1&1&1&0&0&0&0&0&0&0\\
1&1&1&1&0&0&0&0&1&1&1&1&0&0&0\\
1&1&1&0&1&0&0&0&1&0&0&0&1&1&1\\
1&0&0&1&0&1&1&0&1&1&0&0&1&1&0\\
0&1&0&0&1&1&0&1&1&1&1&0&1&0&0\\
1&0&0&0&1&0&1&1&1&0&1&1&0&1&0\\
1&1&0&0&0&1&1&0&0&0&1&1&1&0&1\\
0&0&1&1&0&0&1&1&1&0&1&0&1&0&1\\
1&0&0&1&1&0&0&1&0&1&0&1&1&0&1\\
0&0&1&0&1&1&1&0&1&1&0&1&0&0&1\\
0&1&0&1&1&0&1&0&0&1&1&0&0&1&1\\
0&0&1&1&1&1&0&0&0&0&1&1&1&1&0\\
0&1&1&0&0&0&1&1&0&1&0&1&1&1&0\\
0&1&0&1&0&1&0&1&1&0&0&1&0&1&1\\
1&0&1&0&0&1&0&1&0&1&1&0&0&1&1
\endVmatrix
\mytag{3.1}
$$
The first row of the matrix \mythetag{3.1} is partitioned
into two groups --- the group of ones and the group of zeros next to it:
$$
\hskip -2em
\aligned
&\hphantom{a}\overbrace{\hphantom{aaaaaaaaaaaaaaaaaaaaa}}^{8}
\hphantom{\,a}\overbrace{\hphantom{aaaaaaaaaaaaaaaaaa}}^{7}
\\
\vspace{-1ex}
&\Vmatrix
1&1&1&1&1&1&1&1&0&0&0&0&0&0&0\\
\vdots & \vdots & \vdots & \vdots & \vdots & \vdots & \vdots & \vdots & \vdots & \vdots 
& \vdots & \vdots & \vdots & \vdots & \vdots
\endVmatrix
\endaligned
\mytag{3.2}
$$
We use the Maxima programming language (see \mycite{10}) in order to present the partition structure \mythetag{3.2}. Each of the two groups \mythetag{3.2} is presented by two numbers. 
The first number is the ordering number of the group in the row. The second number is the number of elements in the group. As a result the first row of the matrix \mythetag{3.1}
is presented through the following code:\par
\medskip
\parshape 1 10pt 350pt 
{\tt\noindent\darkred{r1:[[0,8],[1,7]]\$}}\par
\medskip
\noindent
Square brackets in Maxima programming language delimit lists. Therefore each row of the 
matrix \mythetag{3.1} is presented as a list of lists.\par 
     In the second row of the matrix \mythetag{3.1} we see the following four groups:
$$
\hskip -2em
\aligned
&\hphantom{a}\overbrace{\hphantom{aaaaaaaaaa}}^{4}
\hphantom{a\!}\overbrace{\hphantom{aaaaaaaaaa}}^{4}
\hphantom{a\!}\overbrace{\hphantom{aaaaaaaaaa}}^{4}
\hphantom{a}\overbrace{\hphantom{aaaaaaa}}^{3}
\\
\vspace{-1ex}
&\Vmatrix
1&1&1&1&0&0&0&0&1&1&1&1&0&0&0\\
\vdots & \vdots & \vdots & \vdots & \vdots & \vdots & \vdots & \vdots & \vdots & \vdots 
& \vdots & \vdots & \vdots & \vdots & \vdots
\endVmatrix
\endaligned
\mytag{3.3}
$$
The first two groups in \mythetag{3.3} are subordinate to the first group in \mythetag{3.1},
the last two groups of \mythetag{3.3}  are subordinate to the second group in \mythetag{3.1}.
Groups with even numbers correspond to ones in the matrix and the groups with odd numbers
correspond to zeros. As a result we have the following code for the second row of the matrix
\mythetag{3.1}:
\medskip
\parshape 1 10pt 350pt 
{\tt\noindent\darkred{r2:[[0,4],[1,4],[2,4],[3,3]]\$}}\par
\medskip
\noindent
Upon defining each row, the matrix $H$ from \mythetag{3.1} is presented as the list of
its rows:
\medskip
\parshape 1 10pt 350pt 
{\tt\noindent\darkred{H:[r1,r2,r3,r4,r5,r6,r7,r8,r9,r10,r11,r12,r13,r14,r15]\$}}\par
\medskip
     The upper estimate for the number of groups in the $n$-th row of a matrix is $2^n$. 
However the actual number of groups in the row description is much smaller since the 
groups with zero elements are not explicitly listed, though they implicitly assumed. 
Here is the explicit group lists presentation for the matrix \mythetag{3.1}:
\medskip
\parshape 1 10pt 350pt 
{\tt\noindent\darkred{H:[[[0,8],[1,7]],}\newline
\darkred{\ \ \ [[0,4],[1,4],[2,4],[3,3]],}\newline
\darkred{\ \ \ [[0,3],[1,1],[2,1],[3,3],[4,1],[5,3],[6,3]],}\newline
\darkred{\ \ \ [[0,1],[1,2],[2,1],[5,1],[6,2],[7,1],[8,1],[10,1],[11,2],}\newline
\darkred{\ \ \ \ [12,2],[13,1]],}\newline
\darkred{\ \ \ [[1,1],[2,1],[3,1],[5,1],[10,1],[12,1],[13,1],[14,1],[16,1],}\newline
\darkred{\ \ \ \ [20,1],[22,1],[23,1],[24,1],[25,1],[27,1]],}\newline
\darkred{\ \ \ [[2,1],[5,1],[7,1],[11,1],[20,1],[25,1],[26,1],[28,1],[32,1],}\newline
\darkred{\ \ \ \ [41,1],[44,1],[46,1],[49,1],[50,1],[55,1]],}\newline
\darkred{\ \ \ [[4,1],[10,1],[15,1],[23,1],[41,1],[50,1],[52,1],[57,1],[65,1],}\newline
\darkred{\ \ \ \ [83,1],[88,1],[92,1],[98,1],[101,1],[110,1]],}\newline
\darkred{\ \ \ [[9,1],[21,1],[30,1],[46,1],[83,1],[101,1],[104,1],[114,1],}\newline
\darkred{\ \ \ \ [130,1],[167,1],[176,1],[185,1],[196,1],[203,1],[220,1]],}\newline
\darkred{\ \ \ [[18,1],[43,1],[61,1],[92,1],[166,1],[203,1],[209,1],[228,1],}\newline
\darkred{\ \ \ \ [261,1],[334,1],[353,1],[370,1],[392,1],[407,1],[440,1]],}\newline
\darkred{\ \ \ [[37,1],[87,1],[122,1],[185,1],[332,1],[406,1],[418,1],[457,1],}\newline
\darkred{\ \ \ \ [522,1],[668,1],[707,1],[740,1],[785,1],[815,1],[880,1]],}\newline
\darkred{\ \ \ [[75,1],[174,1],[245,1],[370,1],[664,1],[813,1],[836,1],}\newline
\darkred{\ \ \ \ [915,1],[1045,1],[1336,1],[1414,1],[1481,1],[1571,1],}\newline
\darkred{\ \ \ \ [1630,1],[1760,1]],}\newline
\darkred{\ \ \ [[151,1],[349,1],[490,1],[740,1],[1328,1],[1626,1],[1673,1],}\newline
\darkred{\ \ \ \ [1831,1],[2091,1],[2673,1],[2828,1],[2962,1],[3142,1],}\newline
\darkred{\ \ \ \ [3260,1],[3521,1]],}\newline
\darkred{\ \ \ [[303,1],[698,1],[980,1],[1481,1],[2657,1],[3253,1],[3346,1],}\newline
\darkred{\ \ \ \ [3662,1],[4183,1],[5346,1],[5657,1],[5924,1],[6284,1],}\newline
\darkred{\ \ \ \ [6520,1],[7043,1]],}\newline
\darkred{\ \ \ [[607,1],[1396,1],[1961,1],[2962,1],[5315,1],[6506,1],[6693,1],}\newline
\darkred{\ \ \ \ [7324,1],[8366,1],[10693,1],[11315,1],[11848,1],[12569,1],}\newline
\darkred{\ \ \ \ [13040,1],[14086,1]],}\newline
\darkred{\ \ \ [[1214,1],[2793,1],[3922,1],[5925,1],[10631,1],[13012,1],}\newline
\darkred{\ \ \ \ [13387,1],[14648,1],[16733,1],[21386,1],[22630,1],[23697,1],}\newline
\darkred{\ \ \ \ [25139,1],[26080,1],[28172,1]]]\$}
}\par
\medskip
The above presentation of the matrix \mythetag{3.1} looks more complicated than the
regular presentation of matrices in Maxima. However this presentation is more 
convenient from the algorithmic point of view. Below we shall call it the ``group 
lists presentation'' or the ``partition lists presentation. 
\head
4. An algorithm for generating Hadamard matrices. 
\endhead
     The partition lists presentation of Hadamard matrices is the basis for the 
algorithm suggested below. We shall not describe the algorithm in full details
verbally. Instead we provide the source code of it using Maxima programming language
(see \mycite{10}). Hadamard matrices are generated by the following code:
\medskip
\parshape 1 10pt 350pt 
{\tt\noindent\darkred{HM\_size:m\$}\newline
\darkred{HM\_quarter:(HM\_size+1)/4\$}\newline
\darkred{q:HM\_quarter\$}\newline
\darkred{HM\_row[1]:[[0,2*q],[1,2*q-1]]\$}\newline
\darkred{HM\_row[2]:[[0,q],[1,q],[2,q],[3,q-1]]\$}\newline
\darkred{HM\_matrix\_num:1\$}\newline
\darkred{HM\_stream:openw("output\_file.txt")\$}\newline
\darkred{HM\_make\_row(3)\$}\newline
\darkred{close(HM\_stream)\$}
}
\medskip
\noindent
Most of this code initializes global variables. The first variable {\tt\darkred{HM\_size}}
defines the size of Hadamard matrices to be generated. It is given by $m$ which is any
positive integer greater than or equal to $3$ and obeying the condition
$$
m=3\kern -6pt \mod 4.
$$
Otherwise an error message is generated.\par
    The partition lists of the first two rows of the matrix are predefined. They are
stored in the variables {\tt\darkred{HM\_row[1]}} and {\tt\darkred{HM\_row[2]}}
(compare with \mythetag{3.2} and \mythetag{3.3}). The third row and all other rows of 
Hadamard matrices are produced by calling the recursive function 
{\tt\darkred{HM\_make\_row()}} with the argument $3$. The whole job is practically done 
by this function. Below is its code:\par 
\medskip
\parshape 1 10pt 350pt 
{\tt\noindent\darkred{HM\_make\_row(i):=block}\newline
\darkred{\ ([n,s,k,l,q,dummy,kk,y,dpnd,indp,nrd,nri,r,kr,qq,eq,eq\_list,j,}\newline
\darkred{\ \ LLL,RLL,RVV,RRV,subst\_list],}\newline
\darkred{\ \ if not integerp(HM\_size) or HM\_size<3 or mod(HM\_size,4)\#3}\newline
\darkred{\ \ \ then}\newline
\darkred{\ \ \ \ (}\newline
\darkred{\ \ \ \ \ print(printf(false,"Error:\ m=\~{}a is incorrect size for }\newline
\darkred{\ \ \ \ \ Hadamard matrices",HM\_size)),}\newline
\darkred{\ \ \ \ \ return(false)}\newline
\darkred{\ \ \ \ ),}\newline
\darkred{\ \ if HM\_size=3}\newline
\darkred{\ \ \ then}\newline
\darkred{\ \ \ \ (}\newline
\darkred{\ \ \ \ \ HM\_row[2]:[[0,1],[1,1],[2,1]],}\newline
\darkred{\ \ \ \ \ HM\_row[3]:[[1,1],[2,1],[4,1]],}\newline
\darkred{\ \ \ \ \ HM\_output\_matrix(),}\newline
\darkred{\ \ \ \ \ return(false)}\newline
\darkred{\ \ \ \ ),}\newline
\darkred{\ \ print(printf(false,"i=\~{}a",i)),}\newline
\darkred{\ \ n:length(HM\_row[i-1]),}\newline
\darkred{\ \ k:makelist(100,y,n),}\newline
\darkred{\ \ l:makelist(100,y,n),}\newline
\darkred{\ \ q:makelist(100,y,n),}\newline
\darkred{\ \ for s:1 step 1 thru n do}\newline
\darkred{\ \ \ (}\newline
\darkred{\ \ \ \ l[s]:HM\_row[i-1][s][1],}\newline
\darkred{\ \ \ \ q[s]:HM\_row[i-1][s][2]}\newline
\darkred{\ \ \ ),}\newline
\darkblue{\ \ /*-- prepare the equation list --*/}\newline
\darkred{\ \ eq\_list:[],}\newline
\darkred{\ \ var\_list:[],}\newline
\darkred{\ \ eq:0,}\newline
\darkred{\ \ for s:1 step 1 thru n do }\newline
\darkred{\ \ \ (}\newline
\darkred{\ \ \ \ eq:eq+HM\_V[i][s],}\newline
\darkred{\ \ \ \ var\_list:endcons(HM\_V[i][s],var\_list)}\newline
\darkred{\ \ \ ),}\newline
\darkred{\ \ eq\_list:endcons(eq=2*HM\_quarter,eq\_list),}\newline
\darkred{\ \ qq:1,}\newline
\darkred{\ \ for j:i-1 step -1 thru 1 do}\newline
\darkred{\ \ \ (}\newline
\darkred{\ \ \ \ eq:0,}\newline
\darkred{\ \ \ \ for s:1 step 1 thru n do }\newline
\darkred{\ \ \ \ \ if evenp(floor(l[s]/qq)) then eq:eq+HM\_V[i][s],}\newline
\darkred{\ \ \ \ eq\_list:endcons(eq=HM\_quarter,eq\_list),}\newline
\darkred{\ \ \ \ qq:qq*2}\newline
\darkred{\ \ \ ),}
}\par
\parshape 1 10pt 350pt 
{\tt\noindent\darkblue{\ \ /*----- solve the equations -----*/}\newline
\darkred{\ \ linsolve\_params:false,}\newline
\darkred{\ \ HM\_SOL[i]:linsolve(eq\_list,var\_list),}\newline
\darkred{\ \ LLL:map(lhs,HM\_SOL[i]),}\newline
\darkred{\ \ dpnd:[],}\newline
\darkred{\ \ for r:1 step 1 thru length(LLL) do }\newline
\darkred{\ \ \ dpnd:endcons(args(LLL[r])[1],dpnd),}\newline
\darkred{\ \ RLL:map(rhs,HM\_SOL[i]),}\newline
\darkred{\ \ RVV:map(listofvars,RLL),}\newline
\darkred{\ \ RRV:\{\},}\newline
\darkred{\ \ for r:1 step 1 thru length(RVV) do }\newline
\darkred{\ \ \ RRV:union(RRV,setify(RVV[r])),}\newline
\darkred{\ \ RRV:listify(RRV),}\newline
\darkred{\ \ indp:[],}\newline
\darkred{\ \ for r:1 step 1 thru length(RRV) do }\newline
\darkred{\ \ \ indp:endcons(args(RRV[r])[1],indp),}\newline
\darkblue{\ \ /*-- initiate the multiindex loop --*/}\newline
\darkred{\ \ nrd:length(dpnd),}\newline
\darkred{\ \ nri:length(indp),}\newline
\darkred{\ \ kr::makelist(0,y,nri+1),}\newline
\darkred{\ \ for dummy:1 step 1 while kr[nri+1]=0 do}\newline
\darkred{\ \ \ (}\newline
\darkred{\ \ \ \ subst\_list:[],}\newline
\darkred{\ \ \ \ for r:1 step 1 thru nri do}\newline
\darkred{\ \ \ \ \ (}\newline
\darkred{\ \ \ \ \ \ s:indp[r],}\newline
\darkred{\ \ \ \ \ \ k[s]:kr[r],}\newline
\darkred{\ \ \ \ \ \ kk[s]:q[s]-k[s],}\newline
\darkred{\ \ \ \ \ \ subst\_list:endcons(HM\_V[i][s]=k[s],subst\_list)}\newline
\darkred{\ \ \ \ \ ),}\newline
\darkred{\ \ \ \ for r:1 step 1 thru nrd do}\newline
\darkred{\ \ \ \ \ (}\newline
\darkred{\ \ \ \ \ \ s:dpnd[r],}\newline
\darkred{\ \ \ \ \ \ k[s]:psubst(subst\_list,RLL[r]),}\newline
\darkred{\ \ \ \ \ \ kk[s]:q[s]-k[s]}\newline
\darkred{\ \ \ \ \ ),}\newline
\darkblue{\ \ \ \ /*----- create a new row -----*/}\newline
\darkred{\ \ \ \ HM\_row[i]:[],}\newline
\darkred{\ \ \ \ for s:1 step 1 thru n do}\newline
\darkred{\ \ \ \ \ (}\newline
\darkred{\ \ \ \ \ \ if k[s]\#0 then HM\_row[i]:endcons([2*l[s],k[s]],HM\_row[i]),}\newline
\darkred{\ \ \ \ \ \ if kk[s]\#0 then HM\_row[i]:endcons([2*l[s]+1,kk[s]],HM\_row[i])}\newline
\darkred{\ \ \ \ \ ),}\newline
\darkred{\ \ \ \ if HM\_sc\_prods\_ok(i)}\newline
\darkred{\ \ \ \ \ then}\newline
\darkred{\ \ \ \ \ \ if i=n}\newline
\darkred{\ \ \ \ \ \ \ then HM\_output\_matrix()}\newline
\darkred{\ \ \ \ \ \ \ else HM\_make\_row(i+1),}\newline
\darkblue{\ \ \ \ /*--- increment the multiindex ---*/}\newline
\darkred{\ \ \ \ if nri=0 then kr[nri+1]:1, }\newline
\darkred{\ \ \ \ for r:1 step 1 thru nri do}\newline
\darkred{\ \ \ \ \ (}\newline
\darkred{\ \ \ \ \ \ if r=1 then kr[1]:kr[1]+1,}\newline
\darkred{\ \ \ \ \ \ s:indp[r],}\newline
\darkred{\ \ \ \ \ \ if kr[r]>q[s] then (kr[r]:0, kr[r+1]:kr[r+1]+1)}\newline
\darkred{\ \ \ \ \ )}\newline
\darkred{\ \ \ )}\newline
\darkred{\   )\$}
}\par
\medskip
     There are two auxiliary functions which are called from within the above code. 
One of them is used in order to output generated matrices.
\medskip
\parshape 1 10pt 350pt 
{\tt\noindent\darkred{HM\_output\_matrix():=block}\newline
\darkred{\ ([s,LL],}\newline
\darkred{\ \ LL:[],}\newline
\darkred{\ \ for s:1 step 1 thru HM\_size do LL:endcons(HM\_row[s],LL),}\newline
\darkred{\ \ printf(HM\_stream,%
"HM\_\~{}a\_\~{}a:\~{}a\$\~{}\%",HM\_size,HM\_matrix\_num,LL),}\newline
\darkred{\ \ HM\_matrix\_num:HM\_matrix\_num+1}\newline
\darkred{\ )\$}
}\par
\medskip
\noindent 
The second function verifies if the row data prepared for output are correct.
\medskip
\parshape 1 10pt 350pt 
{\tt\noindent\darkred{HM\_sc\_prods\_ok(i):=block}\newline
\darkred{\ ([n,result,ss,s,j,qq],}\newline
\darkred{\ \ n:length(HM\_row[i]),}\newline
\darkred{\ \ result:true,}\newline
\darkred{\ \ ss:0,}\newline
\darkred{\ \ for s:1 step 1 thru n do}\newline
\darkred{\ \ \ (}\newline
\darkred{\ \ \ \ result:result and (HM\_row[i][s][2]>0)}\newline
\darkred{\ \ \ ),}\newline
\darkred{\ \   return(result)}\newline
\darkred{ )\$}
}\par
\medskip
\head
5. Running the code and results. 
\endhead
The above code was run in Maxima, version 5.42.2, on Linux platform of
Ubuntu 16.04 LTS using laptop computer DEXP Atlas H161 with the processor
unit Intel Core i7-4710MQ. Here we focus on performance of the algorithm.\par
    {\bf The case $m=3$} is trivial. In this case the algorithm terminated 
instantly and produced exactly one Hadamard matrix.\par
    {\bf The case $m=7$} is less trivial. In this case the algorithm also
terminated instantly and produced 25 matrices. All of them were tested and
turned out to be correct $7\times 7$ Hadamard matrices in $\{0,1\}$ 
presentation.\par
    {\bf The case $m=11$}. In this case the algorithm terminated upon running
for 3 minutes and 45 seconds. It produced 60481 matrices. Ten of these matrices
randomly chosen were tested and turned out to be correct $11\times 11$ Hadamard 
matrices in $\{0,1\}$ presentation. The matrix production rate is
$$
v=16128\text{\ matrices/minute.}
$$\par
    {\bf The case $m=15$}. In this case the algorithm did not terminate during
observably short time. Upon running for 16 minutes and 49 seconds it produced
162500 matrices. Ten of these matrices randomly chosen were tested and 
turned out to be correct $15\times 15$ Hadamard matrices in $\{0,1\}$ presentation.
The production rate is
$$
v=9663\text{\ matrices/minute.}
$$\par
    {\bf The case $m=19$}. Again the algorithm did not terminate during observably short 
time. Upon running for 1 minute and 40 seconds it produced 10000 matrices. Ten 
randomly chosen matrices were tested and passed the test. The rate is
$$
v=6000\text{\ matrices/minute.}
$$\par
    {\bf The case $m=23$}  is similar to the previous one. The algorithm did not 
terminate during observably short time. 10000 matrices were generated for 2 minutes
and 43 seconds. Ten randomly chosen matrices were tested. They turned out to be correct
$23\times 23$ Hadamard matrices in $\{0,1\}$ presentation. 
The production rate is
$$
v=3680\text{\ matrices/minute.}
$$\par
    {\bf The case $m=27$} is absolutely different. In this case the algorithm ran
overnight for several hours but did not produce any matrices. So $m=27$ is a practical
limit for the algorithm.\par
     One of the features of the present algorithm is that it solves linear equations 
arising from the form of the Gram matrix \mythetag{2.6} using Maxima's 
{\tt\darkred{linsolve}} function instead of scanning over the ranges of variables. 
However it does not solve inequalities in this manner 
(see {\tt\darkred{HM\_row[i][s][2]>0}} in the code of the function {\tt\darkred{HM\_sc\_prods\_ok}} above). Probably using some module for solving linear 
inequalities along with linear equations could improve the algorithm and $m=27$ would
not be a limit for it any more.
\head
4. Dedicatory.
\endhead
     This paper is dedicated to my sister Svetlana Abdulovna Sharipova. 

\enddocument
\end